\documentclass[11pt]{article}
\usepackage{amsmath}
\usepackage{amssymb}
\usepackage{amscd}

\newcommand{\R}{\mathbb{R}}

\newcommand{\abs}[1]{\left\vert #1\right\vert}

\newcommand{\E}{\mathbb{E}}

\newcommand{\eps}{\varepsilon}

\def \EE {\mathcal{E}}

\newcommand{\vc}{{\rm vc}}

\newcommand{\bi}{B\bigl(L_\infty(\Omega)\bigr)}

\newcommand{\vr}[1]{{\rm vr}(#1)}

\newtheorem{Theorem}{Theorem}[section]
\newtheorem{Lemma}[Theorem]{Lemma}

\newtheorem{Definition}[Theorem]{Definition}

\newtheorem{Proposition}[Theorem]{Proposition}
\newtheorem{Corollary}[Theorem]{Corollary}

\numberwithin{equation}{section} 

\def \proof {\noindent {\bf Proof.}\ \ }

\def \endproof
{{\mbox{}\nolinebreak\hfill\rule{2mm}{2mm}\par\medbreak}}

\def \R {\mathbb{R}}

\def \E {\mathbb{E}}
\def \P {\mathbb{P}}

\def \e {\varepsilon}
\def \eps {\varepsilon}
\def \d {\delta}

\def \l {\lambda}

\def \s {\sigma}

\def \< {\langle}
\def \> {\rangle}

\def \ptwo {{\psi_2}}
\def \ptwon {{\psi_2^n}}

\def \vr {{\rm vr}}

\def \bi {B\bigl(L_\infty(\Omega)\bigr)}
\def \vc {{\rm vc}}
\def \dd {{\tilde{\delta}}}
\begin{document}
\title {Remarks on the Geometry of Coordinate Projections in $\R^n$}
\author {S. Mendelson\footnote{Research School of
Information Sciences and Engineering, The Australian National
University, Canberra, ACT 0200, Australia, e-mail:
shahar.mendelson@anu.edu.au} \and
   R. Vershynin\footnote{
   Department of Mathematical Sciences,
   University of Alberta,
   Edmonton, Alberta T6G 2G1, Canada,
   e-mail: rvershynin@math.ualberta.ca}}
\date{}

\maketitle
\begin{abstract}
We study geometric properties of coordinate projections. Among
other results, we show that if a body $K \subset \R^n$ has an
``almost extremal" volume ratio, then it has a projection of
proportional dimension which is close to the cube. We also
establish a sharp estimate on the shattering dimension of the convex
hull of a class of functions in terms of the shattering dimension
of the class itself.
\end{abstract}

\section{Introduction}
In this article we present several results on coordinate
projections. The majority of this article is devoted to new
applications of the entropy inequality established in \cite{MenVer},
which, roughly speaking, states that if a set of functions has a
large entropy in $L_2$, it must have a coordinate projection
which contains a large cube.

\begin{Definition} \label{def:fat}
We say that a subset $\s$ of $\Omega$ is {\bf $t$-shattered} by a class
of real-valued functions $F$ 
if there exists a level function $h$ on $\s$ such that, given
any subset $\s'$ of $\s$, one can find a function $f \in F$ with
$f(x)  \le  h(x) - t$ if $x \in \s'$ and $f(x)  \ge  h(x) + t$ if $x
\in \s \setminus \s'$.

The {\bf shattering dimension} of $A$, denoted by
$\vc(F,\Omega,t)$ after Vapnik and Chervonenkis, is the maximal
cardinality of a subset of $\Omega$ which is $t$-shattered by $F$.
In cases where the underlying space is clear we denote the
shattering dimension by $\vc(F,t)$.
\end{Definition}

\begin{Theorem} \label{thm:main}
 Let $A$ be a class of functions bounded by $1$,
  defined on a set $\Omega$.
  Then for every probability measure $\mu$ on $\Omega$,
  \begin{equation}                          \label{combinatorial}
    N(F, t, L_2(\mu))
    \le  \Big( \frac{2}{t} \Big)^{K \cdot \vc(F,\, c t)},
     \ \ \ \ 0 < t < 1,
  \end{equation}
  where $K$ and $c$ are positive absolute constants.
\end{Theorem}
Every $F \subset \R^n$ can be identified with a class of functions
on $\{1,...,n\}$ in the natural way: $v(i)=v_i$ for $v \in F$. If
we take $\mu$ to be the probability counting measure on
$\{1,...,n\}$ then \eqref{combinatorial} states that for any
$0<t<1$,
\begin{equation*}
N(F,t\sqrt{n}B_2^n) \leq \Big( \frac{2}{t} \Big)^{K \cdot \vc(F,\,
c t)}.
\end{equation*}
Note that if $F$ happens to be convex and symmetric with respect
to the origin, then $\vc(F,t)$ is the maximal cardinality of a
subset $\s$ of $\{1, \ldots, n\}$ such that $P_\s(F) \supset
[-t, t]^\s$.

We apply this result to study convex bodies whose volume ratio is
almost maximal. Recall that the volume ratio, introduced by Szarek
and Tomczak-Jaegermann \cite{S,ST}, is defined as $\vr(D) =
(|D|/|\EE|)^{1/n}$, where $\EE$ is the ellipsoid of maximal volume
contained in $D$.

The minimal volume ratio of a symmetric convex body in $\R^n$ is
$1$ and is attained by the Euclidean ball; the maximal is of the
order of $\sqrt{n}$ and is attained by the cube $B_\infty^n$
\cite{B}. This pair of extremal bodies is unique up to a linear
transformation. Indeed, the uniqueness of the minimizer is
immediate, while the fact that cube is the unique maximizer was
established in \cite{Ba}.

The isomorphic version of this fact -- describing the bodies whose
volume ratio is of order either $1$ or $\sqrt{n}$ is of particular
interest. The question is whether such bodies inherit any
structure from the Euclidean ball or, respectively, from the cube.

If the volume ratio of a body $D$ in $\R^n$ is bounded by a constant,
then by the Volume Ratio Theorem \cite{ST}, $D$ has a section of
dimension proportional to $n$, which is well isomorphic to the
Euclidean ball.

On the other hand, if $\vr(D)$ is of order of $\sqrt{n}$, then by
\cite{R} and \cite{V}, $D$ has a section of dimension proportional
to $\sqrt{n}$, which is well isomorphic to the cube, and the order
of $\sqrt{n}$ in the dimension can not be improved
(the dual of Gluskin's polytope is such an example -- see section \ref{evr}).
However, we show in section \ref{evr} that there exists a {\em
projection} of $D$ of dimension proportional to $n$, which is well
isomorphic to the cube.

Two other applications we present are based on the following
corollary of Theorem \ref{thm:main}.

\begin{Theorem} \label{thm:Talagrand} \cite{MenVer}
  Let $F$ be a class of functions bounded by $1$,
  defined on a finite set $I$ of cardinality $n$.
  Then the gaussian process indexed by $F$, $X_f=\sum_{i=1}^n g_i f(i)$
  satisfies
  $$
  E=\E \sup_{f \in F} X_f \le  K \sqrt{n}
       \int_{cE/n}^{1} \sqrt{\vc(F,t) \cdot \log (2/t)}\; dt,
  $$
  where $K$ and $c$ are absolute constants.
\end{Theorem}

One application we present is a comparison of the average $\E \|
\sum_{i=1}^n \e_i x_i\|$ to the minimum over all choices of signs,
$\min \| \sum_{i=1}^n \pm x_i\|$. As a consequence, we compare the
type $2$ and the infratype $2$ constants of a Banach space.

Then, we establish a sharp estimate on the shattering dimension of
a convex hull of a class of functions, based on the shattering
dimension of the class itself. Namely, we show that for every
$\eps>0$,

$$
\vc({\rm conv}(F),\eps) \le (C / \e)^2 \cdot \vc(F,c\eps),
$$
where $c$ and $C$ are absolute constants.

The final question we address is when a random coordinate
projection an ``almost isometry". Let $(\Omega,\mu)$ be a
probability space and $f \in L_2(\mu)$, and for simplicity, assume
that $\Omega=\{1,...,n\}$ and that $\mu$ is the uniform
probability measure on $\Omega$. For every $\eps>0$, our aim is to
find ``many" sets $\sigma \subset \{1,...,n\}$ of small
cardinality such that the natural coordinate projection $P_\s f $
satisfies
\begin{equation} \label{eq:isometry}
(1-\eps)\|f\|_{L_2^n} \leq \|P_\s f\|_{L_2^\s} \leq
(1+\eps)\|f\|_{L_2^n},
\end{equation}
where $L_2^k$ is the $L_2$ space defined on $\{1,...,k\}$ with
respect to the uniform probability measure.

By a standard concentration argument, if $\|f\|_\infty \leq 1$,
then with high probability a random coordinate projection of
dimension $C/\eps^2$ is an almost isometry in the sense of
\eqref{eq:isometry}. We will show that the uniform boundedness of
$f$ can be relaxed; it suffices to assume that $\|f\|_{\psi_2}
\leq 1$, where $\| \ \|_{\psi_2}$ is the Orlicz norm generated by
the function $e^{t^2}-1$. In this case, a random coordinate projection of
dimension $C/\eps^2$ will be an almost isometry as in
\eqref{eq:isometry} with high probability. Although this result is
relatively easy, we decided to present it because it gives hope
that the conditions in stronger concentration inequalities (e.g.
Talagrand's concentration inequality for empirical processes
\cite{T 94,L}) can also be relaxed. As an application, we obtain a
coordinate version of the Johnson-Lindenstrauss ``Flattening"
Lemma \cite{JL}.

Finally, we turn to some notational conventions. Throughout, all
absolute constants are denoted by $c$, $C$, $k$ and $K$. Their
values may change from line to line or even within the same line.
We denote $a \sim b$ if there are absolute constants $c$ and $C$
such that $cb \leq a \leq Cb$.

\qquad

{\bf ACKNOWLEDGEMENTS:} The first author acknowledges partial
support by an Australian Research Council Discovery Grant. The
second author thanks Nicole Tomczak-Jaegermann for her constant
support. A part of this work was done when the second author was
visiting the Research School of Information Sciences and
Engineering at The Australian National University, which he thanks
for its hospitality. He also acknowledges a support from the
Pacific Institute of Mathematical Sciences and is grateful to the
Department of Mathematical Sciences of the University of Alberta
for its hospitality.

\section{Extremal volume ratios}    \label{evr}

The volume ratio of a convex body $D$ in $\R^n$ is
defined as
$$
\vr(D)  =  \inf \Big( \frac{|D|}{|\EE|} \Big)^{1/n},
$$
where $| \ |$ denotes the volume in $\R^n$, and the infimum is
over all ellipsoids $\EE$ contained in $D$. This important
invariant was introduced by Szarek and Tomczak-Jaegermann (see
\cite{S}, \cite{ST} or \cite{P}).

The bodies with extremal volume ratios are the Euclidean ball and
the cube -- and these are the only extreme bodies up to a linear
transformation (for the uniqueness of the cube, see \cite{Ba}).
One can show that for every convex symmetric body in $\R^n$,
\begin{equation}                                \label{range vr}
1  =  \vr(B_2^n)  \le  \vr(D)  \le  \vr(B_\infty^n),
\end{equation}
(see \cite{B}), while direct computation shows that
$\vr(B_\infty^n) \leq C\sqrt{n}$ and the best value of the
constant is $C = 2 / \sqrt{\pi e}$.

Often, one encounters bodies whose volume ratio is {\em almost}
extremal, i.e. close to one of sides of \eqref{range vr}. The
problem is whether such a body inherits properties of the extremal
bodies, the Euclidean ball or the cube.

If $\vr(D) \le A$, then by the Volume Ratio Theorem \cite{ST}, $D$
has a section of dimension $k = n/2$ which is $c A^2$-isomorphic
to the Euclidean ball $B_2^k$, and this result is asymptotically
sharp.

On the opposite side of the scale, if $\vr(D) \ge A^{-1}
\sqrt{n}$, $D$ has a section of $D$ of dimension $k = c(A)
\sqrt{n}$ which is $C(A) \log n$-isomorphic to the cube
$B_\infty^k$ \cite{R,V}. It is not known whether the logarithmic
term can be eliminated, but the order of $\sqrt{n}$ in the
dimension is optimal, as was noticed in \cite{GTT}.
Indeed, by an argument of Figiel and Johnson (see
\cite{FJ}, cor. 3.2), a random subspace $E \subset \ell_\infty^n$
(and thus a dual of Gluskin's space) of dimension at least
$n/2$ satisfies that for any $F \subset E$, $gl(F) \geq c{\rm
dim}(F)/\sqrt{n}$, where $gl(F)$ is the Gordon-Lewis constant of
$F$, and $c$ is a suitable absolute constant. By \cite{GL}, $gl(F)
\leq {\rm unc}(F)$, where ${\rm unc}(F)$ is the least
unconditionality constant of a basis of $F$. Since ${\rm
unc}(\ell_\infty^k)=1$ then
$$
d(F,\ell_\infty^k) \geq {\rm unc}(F) \geq \frac{c{\rm
dim}(F)}{\sqrt{n}},
$$
and thus, if $F$ is $2$-isomorphic to $\ell_\infty^k$ then ${\rm
dim}(F) \leq c'\sqrt{n}$.

Our next result shows that $D$ has a {\em projection} of dimension
proportional to $n$ which is $cA$-isomorphic to the cube
$B_\infty^k$.

\begin{Theorem} \label{thm:extremal}
There are absolute constants $C$ and $c$ for which the following
holds. If $D$ is a convex symmetric body in $\R^n$ for which
$\vr(D) \geq A^{-1}\sqrt{n}$, then there exists a projection $P$
of rank $k  \ge  c n / \log A$ such that
$$
d(PK, B_\infty^k)  \le  CA.
$$
\end{Theorem}

To prove the Theorem, recall the notion of the {\it cubic ratio}
\cite{B}. For every ball $D \subset \R^n$ one defines
$$
{\rm cr}(D) = \inf \Big( \frac{|B_\infty^n|}{|TD|}
\Big)^{\frac{1}{n}},
$$
where the infimum is over all linear invertible operators $T$ on
$\R^n$ such that $TD \subset B_\infty^n$.

\begin{Lemma} \label{lemma:ball} \cite{B}
There are absolute constants $c$ and $C$ such that for every
integer $n$ and every convex symmetric body $D \subset \R^n$,
\begin{equation*}
c\sqrt{n} \leq {\rm vr}(D) \cdot {\rm cr}(D) \leq C\sqrt{n}.
\end{equation*}
\end{Lemma}

\noindent {\bf Proof of Theorem \ref{thm:extremal}.} Clearly, we
can assume $n$ to be larger than a suitable absolute constant $N$,
which ensures that for every $D \subset \R^n$, $\vr(D) \leq
0.8\sqrt{n}$. Since $\vr(D) \geq A^{-1} \sqrt{n}$, then by Lemma
\ref{lemma:ball}, ${\rm cr}(D) \le CA$. Hence, there is some $T
\in GL_n$ such that
$$
TD \subset B_\infty^n \ \ \text{ and } \ \ |TD|^{\frac{1}{n}} \geq
c/A.
$$
Recall that $c_1^n  \le  |\sqrt{n}B_2^n|  \le  c_2^n$ for some
absolute constants $c_1, c_2$, and thus there exists an absolute
constant $c_3$ such that
$$
2^n = \frac{|TD|}{|c_3 A^{-1}(\sqrt{n}B_2^n)|}.
$$
By a standard volumetric argument, the right-hand side is bounded
by
$$
N(TD, c_3 A^{-1}\sqrt{n}B_2^n),
$$
and by Theorem \ref{thm:main} there are absolute constants $K$ and
$c$ for which
$$
n \le \log N(TD, c_3 A^{-1}\sqrt{n}B_2^n) \leq K \cdot \vc(TD, c
A^{-1}) \log (CA).
$$
Hence, there is a set $\s \subset \{1,...,n\}$, such that $|\s|
\geq n / K \log(CA)$ and
$$
c_1 A^{-1} B_\infty^\s \subset P_{\s} (TD) \subset
B_\infty^\s.
$$
It only remains to note that $\log(CA)  \le  C' \log A$, because
$A \ge 5/4$.
\endproof

\begin{remark}
Since the volume ratio is always greater than $1$, then $A \geq
n^{-1/2}$. Therefore, the dimension of the cubic projection in
Theorem \ref{thm:extremal} is always bounded below by $c n / \log
n$.
\end{remark}

In a very similar way, one can prove the following
\begin{Theorem}
There are absolute constants $C$ and $c$ for which the following
holds. Let $D$ be a convex symmetric body in $\R^n$ for which
$\vr(D) \geq A^{-1}\sqrt{n}$. Then there exists a projection $P$
of rank $k  \ge  cn$ such that
$$
d(PK, B_\infty^k)  \le  C A^2.
$$
\end{Theorem}

\section{Type and Infratype} In this section we improve a result
of M.~Talagrand \cite{T 92} which compares the average over the
$\pm$ signs to the minimum over the $\pm$ signs of $\|\sum_{i=1}^n
\pm x_i\|$.

Recall that a Banach space $X$ has a (gaussian) type $p$ if there
exists some $M > 0$ such that for all $n$ and all sequences of
vectors $(x_i)_{i \le n}$,
\begin{equation}                            \label{Tp}
  \E \Big\| \sum_{i = 1}^n g_i x_i \Big\| \le  M \Big( \sum_{i =
  1}^n \|x_i\|^p \Big)^{\frac{1}{p}}.
\end{equation}
The best possible constant $M$ in this inequality is denoted by
$T_p(X)$. We say that $X$ has infratype $p$ if there exists some
$M
> 0$ such that for all $n$ and all sequences of vectors $(x_i)_{i
\le n}$,
\begin{equation}                            \label{Ip}
  \min_{\eta_i = \pm 1} \Big\| \sum_{i = 1}^n \eta_i x_i \Big\|
  \le  M \Big( \sum_{i = 1}^n \|x_i\|^p \Big)^{\frac{1}{p}}.
\end{equation}
The best possible constant $M$ in this inequality is denoted by
$I_p(X)$.

In \cite{T 92} it was shown that if $1 < p < 2$ then $T_p(X) \le
C_p I_p(X)^2$, where $C_p$ is a constant which depends only on
$p$. It is not known whether the square can be removed. Regarding
the case $p=2$, M. Talagrand recently constructed a symmetric
sequence space which has infratype $2$ but not type $2$ \cite{T
03}. Hence one can not obtain dimension free estimates on
$T_2(X)$ in terms of $I_2(X)$. Our main result in this section is
that is ${\rm dim}(X)=n$ then $T_2(X) \leq CI_2(X) \cdot
\log^{3/2}n$.

We begin with the following fact that allows one to compare
Rademacher and Gaussian averages.

\begin{Lemma}                       \label{comparison pm}
There is an absolute constant $C$ for which the following holds.
Let $x_1, \ldots, x_n$ be vectors in the unit ball
  of a Banach space and let $0<M \leq \sqrt{n}$.
  If $0  < \l  < \log^{-3} (n / M^2)$ and
  $$
  \min_{\eta_i = \pm 1} \Big\| \sum_{i \in \s} \eta_i x_i \Big\|
  \le  M |\s|^{\frac{1}{2}}
$$
for all $\s \subset \{1,...,n\}$ with $|\s|  \le \l n$, then,
  $$
  \E \Big\| \sum_{i=1}^n g_i x_i \Big\|
  \le C M (n / \l)^{1/2}.
  $$
\end{Lemma}

In the proof of Lemma \ref{comparison pm} we require the following
observation from \cite{MS}, that if $\{x_1,...,x_n\} \subset X$ is
$\eps$-shattered by $B_{X^*}$, then for any $(a_i)_{i=1}^n \in
\R^n$,
\begin{equation} \label{eq:cube}
\eps \sum_{i=1}^n |a_i| \leq \Big\|\sum_{i=1}^n a_i x_i \Big\|.
\end{equation}

\noindent{\bf Proof of Lemma \ref{comparison pm}.}\  Clearly we
can assume that the given Banach space is $X = (\R^n, \|\cdot\|)$
and that $(x_i)_{i \le n}$ are the unit coordinate vectors in
$\R^n$. Set $B = B_{X^*}$ and by the hypothesis of the lemma and
\eqref{eq:cube}, $\vc(B, M v^{-1/2}) \le v$ if $0 \le v \le \l n$.
Hence, for any $M (\l n)^{-1/2}  \le t  \le  1$,
\begin{equation}                        \label{VC small}
  \vc(B, t)  \le  (M / t)^2.
\end{equation}
Set
$$
E = \E \Big\|\sum_{i=1}^n g_i x_i \Big\|_X
  = \E \sup_{b \in B} \sum_{i=1}^n g_i \,b(i).
$$
By Theorem \ref{thm:Talagrand}, there are absolute constants $C$
and $c$ such that
$$
E  \le  K \sqrt{n} \int_{c E / n}^1 \sqrt{\vc(B, t) \cdot \log(2 /
t)} \; d t.
$$
If $c E / n  \le  M (\l n)^{-1/2}$, the lemma trivially follows.
Otherwise, if the converse inequality holds, then by \eqref{VC
small} and since $\lambda <1$,
$$
E  \le  K \sqrt{n} \int_{c E / n}^1 (M / t) \sqrt{\log(2 / t)} \;
d t \le  K \sqrt{n} M \cdot \log^{\frac{3}{2}} (n / M^2),
$$
and by the assumption on $\lambda$,
$$
E  \le  K \sqrt{n} M \cdot \log^{\frac{3}{2}} (n / M^2) \le  K
\sqrt{n} M/ \sqrt{\l},
$$
as claimed.
\endproof

Using Lemma \ref{comparison pm}, one can compare the type and
infratype $2$ of a Banach space $X$.

Let $T_2^{(n)} (X)$ and $I_2^{(n)} (X)$ denote the best possible
constants $M$ in \eqref{Tp} and \eqref{Ip} respectively (with $p =
2$). So, $T_2^{(n)} (X)$ and $I_2^{(n)} (X)$ measure the
type/infratype $2$ computed on $n$ vectors. Clearly, $I_2 (X) \le
T_2 (X)$ and $I_2^{(n)} (X)  \le  T_2^{(n)} (X) \leq \sqrt{n}$.

\begin{Theorem}
  Let $X$ be an $n$-dimensional Banach space.
  Then, for every number $0  < \l  < \log^{-3} (n / I_2(X)^2)$,
  $$
  T_2 (X)  \le  C \cdot I_2^{(\l n)} (X)/\sqrt{\l}.
  $$
\end{Theorem}
In particular,
$$
T_2 (X)  \le  I_2 (X)  \cdot  C \log^{\frac{3}{2}} \Big(
\frac{n}{I_2(X)^2} \Big)
  \le  I_2 (X)  \cdot  C \log^{\frac{3}{2}} n.
$$

\proof By \cite{TJ} and \cite{BKT} Theorem 3.1, the Gaussian type
$2$ can be computed on $n$ vectors of norm one. Precisely, this
means that $T_2 (X)$ is the smallest possible constant $M'$ for
which the inequality
$$
\E \Big\| \sum_{i=1}^n g_i x_i \Big\| \le M' n^{1/2}
$$
holds for all vectors $x_1, \ldots, x_n$ of norm one. Now, the
assertion follows from Lemma \ref{comparison pm}.
\endproof

\section{The shattering dimension of convex hulls}
In this section we present a sharp estimate which compares the shattering
dimensions of a class and of its convex hull. To that end, we
connect the shattering dimension to the growth rate of the
expectation of the supremum of the gaussian process $\{X_a,\ a \in
P_\s F\}$ as a function of $|\s|$.

\begin{Definition}
Let $F$ be a class of functions bounded by $1$ and set
$$
\ell_n(F) =\sup_{(x_1,...,x_n) \in \Omega^n} \E_g \sup_{f \in F}
\big| \sum_{i=1}^n g_i f(x_i) \big|,
$$
where $g_1,...,g_n$ are independent, standard gaussian random
variables.
\end{Definition}

Hence, $\ell_n(F)$ is the largest gaussian average associated with
a coordinate projection of $F$ on $n$ points. Since $F \subset
\bi$ then for every $\s=(x_1,...,x_n)$,
$$
P_\s F =\Bigl\{ \bigl(f(x_1),...,f(x_n)\bigr): \ f \in F \Bigr\}
\subset B_\infty^n,
$$
and the largest projection one might encounter is when $P_\s
F=B_\infty^n$, in which case $\ell(P_\s F) \sim n$. We define a
scale-sensitive parameter which measures for every $\eps>0$ the
largest cardinality of a projection which has a ``large"
$\ell$-norm:
$$
t(F,\eps)=\sup \{n: \ \ell_n(F) \geq \eps n \}.
$$

\begin{Theorem} \label{shattering vs. t}
There are absolute constants $K$, $c$ and $c'$ such that for any
$F \subset \bi$ and every $\eps>0$,
\begin{equation*}
\vc(F,c'\eps) \leq t(F,\eps) \leq (K/\eps^2) \cdot \vc(F,c\eps).
\end{equation*}
\end{Theorem}

In the proof, we will use the following wording: the function
$f$ associated to a set $\s'$ in the Definition \ref{def:fat}
will be called the {\bf function that shatters $\s'$}.

\noindent{\bf Proof of Theorem \ref{shattering vs. t}.}\
Assume that
$\{x_1,...,x_n\}$ is $\eps$-shattered by $F$. For every $J \subset
\{x_1,...,x_n\}$, let $f_J$ be the function shattering $J$, and
for each $(\eps_1,...,\eps_n) \in \{-1,1\}^n$ set
$I=\{x_i|\eps_i=1\}$. By the triangle inequality and letting
$f=f_I$, $f'=f_{I^c}$ in the second inequality below,
\begin{align*}
&  \sup_{f \in F} \abs{\sum_{i=1}^n \eps_i f(x_i)} \geq
\frac{1}{2}\sup_{f,f' \in F} \abs{\sum_{i=1}^n
\eps_i\bigl(f(x_i)-f'(x_i)\bigr)}
\\
& \geq \frac{1}{2} \abs{\sum_{i=1}^n
\eps_i\bigl(f_I(x_i)-f_{I^c}(x_i)\bigr)} \geq n\eps.
\end{align*}
Hence,
\begin{equation*}
\sup_{f \in F} \abs{\sum_{i=1}^n \eps_if(x_i)} \geq n \eps,
\end{equation*}
and in particular this holds for the average. The first bound is
evident because of the known connections between gaussian and
Rademacher averages \cite{TJ1}, namely, that there is an absolute
constant $C$ such that for any class $F$ and any set
$\s=(x_1,...,x_n)$,
\begin{align*}
\ell(P_\s F) & =\E_g \| \sum_{i=1}^n g_ie_i \|_{(P_\s F)^\circ}
\geq
C \cdot \E_\eps \|\sum_{i=1}^n \eps_i e_i \|_{(P_\s F)^\circ} \\
& =C \cdot \E_\eps \sup_{f \in F} \big|\sum_{i=1}^n \eps_i f(x_i)
\big|.
\end{align*}
The reverse inequality follows from Theorem \ref{thm:Talagrand} in
a similar way to the proof of Elton's Theorem in \cite{MenVer}. If
$\ell(P_\s F) \geq \eps {n}$, then $\E \sup_{f \in F} X_f \geq
n\eps$, where $X_f=\sum_{i=1}^n g_if(x_i)$. By Theorem
\ref{thm:Talagrand},
$$
n\eps \leq \E\sup_{f \in F} X_f  \leq K \sqrt{n}
       \int_{c\eps}^{1} \sqrt{\vc(P_\s F,t) \cdot \log (2/t)}\; dt.
$$
Set $ v(t)  =  \frac{c_0}{t \log^{1.1} (2 / t)} $ where $c_0 > 0$
is chosen so that $\int_0^{1} v(t) \; dt  =  1$. Hence, there is
some $c\eps \le  t \le  1$ such that $K^2 \cdot \vc(P_\s F,t) \geq
\eps^2n \cdot v^2(t)/\log(2/t)$, implying that
$$
\vc(F,c\eps) \geq \vc(P_\s F, c\eps) \ge \vc(P_\s F, t)  \ge
\frac{c' \eps^2}{t^2 \log^{3.2} (2 / t)} n \ge c{''}n\eps^2.
$$
\endproof
The previous result can be used to estimate the shattering
dimension of a convex hull of a class.

\begin{Corollary} \label{cor:MenVer2}
There are absolute constants $K$ and $c$ such that for any $F
\subset \bi$ and every $\eps>0$,
$$
\vc({\rm conv (F)},\eps) \leq (K/\e)^2 \cdot \vc(F,c\eps).
$$
\end{Corollary}
\proof Since the $\ell$-norm of a set and of its convex hull are
the same, then for any $\eps>0$, $t(F,\eps)=t({\rm
conv}(F),\eps)$. By Theorem \ref{shattering vs. t},
$$
\vc({\rm conv}(F),\eps) \leq t({\rm conv}(F),\eps)=t(F,\eps) \leq
(K / \e)^2 \cdot \vc(F,c\eps).
$$
\endproof
Next, we show that this estimate is sharp, in the sense that the
exponent of $1/\eps^2$ can not be improved. To that end, we
require some properties of the shattering dimension of classes of
linear functionals mentioned before, which was investigated in
\cite{MS}.

If $X$ is a normed space then $B_{X^*}$ can be viewed as a subset
of $L_\infty(B_X)$ in the natural way. It is not difficult to
characterize the shattering dimension in this case.

\begin{Lemma} \label{lemma:shatter}
A set $\{x_1,...,x_n\} \subset B_X$ is $\eps$-shattered by
$B_{X^*}$ if and only if $(x_i)_{i=1}^n$ are linearly independent
and $\eps$-dominate the $\ell_1^n$ unit-vector basis; i.e.,
$$\eps
\sum_{i=1}^n |a_i| \leq \Big\|\sum_{i=1}^n a_i x_i \Big\| \leq \sum_{i=1}^n
|a_i|
$$
for every $a_1,...,a_n \in \R$.

In particular, if $X$ is $n$-dimensional, and if the Banach-Mazur
distance satisfies that $d(X,\ell_1^n) \leq \alpha$, then
$\vc(B_{X^*},B_X,1/\alpha)=n$.
\end{Lemma}

\begin{Corollary} \label{cor:lowerbound}
There exists an absolute constant $k$ for which the following
holds. For every $0<\eps<1/2$ there is a class $F \subset \bi$
such that
\begin{equation*}
\vc({\rm conv}(F),\eps) \geq \frac{k \cdot \vc(F,\eps)}{\eps^2\log
({1}/{\eps})}.
\end{equation*}
\end{Corollary}

\proof For every integer $n$, let $\Omega_n=B_\infty^n$ and set
$F_n=\{e_1,...,e_n\}$, that is, the standard unit vectors in
$\R^n$, when considered as linear functionals on $B_\infty^n$.
Since $|F|=n$, it follows that for every $\eps>0$,
$\vc(F_n,\Omega_n,\eps) \leq \log_2 n$. On the other hand, ${\rm
conv}(F_n)=B_1^n$ when considered as functionals on $B_\infty^n$.
By Lemma \ref{lemma:shatter} applied to $X=\ell_\infty^n$, and
since $d(\ell_\infty^n,\ell_1^n) \leq K\sqrt{n}$ \cite{TJ1}, it is
evident that there is a subset on cardinality $n$ in $B_\infty^n$
which is $k/\sqrt{n}$-shattered by $B_1^n$. Thus, for
$\eps_n=k/\sqrt{n}$,
$$
\vc({\rm conv}(F_n),\Omega_n,\eps_n) \geq \frac{k' \cdot
\vc(F_n,\Omega_n,\eps_n)}{\eps_n^2\log ({1}/{\eps_n})},
$$
from which the proof easily follows.
\endproof

\section{Almost isometric coordinate projections} \label{sec:proj}

Given real-valued function $f$ on a probability space, its
$\psi_p$-norm ($p \ge 1$) is defined as the Orlicz norm
corresponding to the function $\exp(t^p) - 1$. Precisely,
$\|f\|_{\psi_p}$ is the infimum of all numbers $\l$ satisfying $\E
\exp(|f|^p / \l^p)  \le e$. It is possible to compare the $\psi_p$
with other $\psi_q$ norms and the $L_p$ norms. Indeed, one can
show that if $1 \leq p \leq q < \infty$, $\|f\|_{\psi_p} \leq
C_{p,q} \|f\|_{\psi_q}$, and $\|f\|_{L_p} \leq C_p \|f\|_{\psi_1}$
(see, for example, \cite{VW}).

A function $f$ is bounded in the $\ptwo$ norm if and only if $f$
has a subgaussian tail. Namely, if $\|f\|_\ptwo \le 1$ then by
Chebychev's inequality $\P \{ |f|
> t \} \le e^{-t^2 + 1}$ for all $t > 0$. Conversely, if for some $A
\ge 1$ one has $\P \{ |f| > t \} \le A e^{-t^2}$ for all $t
> 1$, then integrating by parts it follows that
$\E \exp(f/2)^2  \le  1 + A/3  \le  2^A$, and by Jensen's
inequality one can conclude that $\|f\|_\ptwo  \le  2 A$ (we did
not attempt here to give the right dependence on $A$).

Another simple but useful fact which follows from Jensen's
inequality is that $\|f\|_\ptwo \le C \, \E \exp (f^2)$, where $C$
is an absolute constant.

We will focus on functions defined on a finite domain, which we
identify with $\{1, \ldots, n\}$, equipped with a uniform measure,
where each atom carries a weight of $1/n$. We denote the $\ptwo$
norm of a function $f$ on this probability space by
$\|f\|_\ptwon$. Since $f$ is defined on $\{1,...,n\}$, we
sometimes identify $f$ with the sequence of scalars
$(f(i))_{i=1}^n$.

We shall use the following standard probabilistic model for random
coordinate projections. Given $0<\delta \leq 1/2$, let
$\delta_1,...,\delta_n$ be selectors, i.e. independent
$\{0,1\}$-valued random variables with mean $\d$. Then $\s = \{i
\;|\; 1 \le i \le n, \  \d_i = 1\}$ is a random subset of the
interval $\{1, \ldots, n\}$ with average cardinality $\d n$.

By Bernstein's inequality \cite{VW}, for every $0<\eps<1$,
$$
\P \Bigl\{ (1-\eps)\|f\|_{L_2^n}^2 \leq \frac{1}{\delta n}
\sum_{i=1}^n \delta_i |f(i)|^2 \leq (1+\eps)\|f\|_{L_2^n}^2
\Bigr\} \geq 1-2\exp\Bigl(-\frac{c\eps^2\delta
n}{\|f\|_\infty}\Bigr),
$$
and by another application of Bernstein's inequality,
\begin{equation} \label{eq:app-Bern}
\P \Bigl \{ \frac{1}{\delta n}\sum_{i=1}^n |\delta_i-\delta| \geq
\eps \Bigr\} \leq 2\exp(-c\eps^2n\delta),
\end{equation}
implying that if $\|f\|_\infty \leq 1$, then with probability at
least $4\exp(-c\eps^2|\s|)$,
$$
(1-\eps)\|f\|_{L_2^n} \leq \|P_\s f\|_{L_2^\s} \leq
(1+\eps)\|f\|_{L_2^n}.
$$

In this section we relax the assumption that $f$ is bounded in the
uniform norm, and assume that $f$ is bounded in the $\ptwo$ norm.

Roughly speaking, we show that for every $1 \leq p<\infty$, the
set of vectors in $S(L_p^n)$ which will be almost isometrically
projected onto $L_p^\s$ are those with a ``small" $\psi_p^n$ norm.

\begin{Proposition} \label{prop}
  Let $(\d_i)_{i=1}^n$ be independent $\{0,1\}$-valued random
  variables with mean $\d > 0$. Set
  $a = (a_i)_{i=1}^n \in \R^n$ and put
  $M = \|a\|_{\psi_1^n}$.
  Then, for every positive number $t < M/2$,
  $$
  \P \Big\{ \sum_{i=1}^n (\d_i - \d) a_i  >  t \d n \Big\}
    \le  \exp \Big( - \frac{c t^2 \d n}{M^2} \Big),
  $$
  where $c$ is an absolute constant.
\end{Proposition}

The proof starts with the following standard lemma.

\begin{Lemma}                               \label{Chebychev}
  Let $Z$ be a random variable and assume that for some $b, \l >
  0$,
  $$
  \E \exp(\l Z)  \le  e^{b^2 \l^2}.
  $$
  Then
  $$
  \P \{ Z > 2 b^2 \l \}  \le  e^{-b^2 \l^2}.
  $$
\end{Lemma}

\proof For $t > 0$,
\begin{align*}
\P \{ Z > t \}
  &= \P \{ \exp ( \l (Z - t) ) > 1 \}
  \le  \E \exp ( \l (Z - t) ) \\
  &= e^{-\l t} \E \exp(\l Z)
  \le  e^{b^2 \l^2 - \l t}.
\end{align*}
Setting $t = 2 b^2 \l$ completes the proof.
\endproof

\noindent {\bf Proof of Proposition \ref{prop}.} By homogeneity,
we can assume that $M=1$, and we shall evaluate $\E \exp
\bigl(t\sum_{i=1}^n (\delta_i-\delta) a_i \bigr)$. To that end,
let $\delta_i'$ be an independent copy of $\delta_i$ and set
$\tilde{\delta}_i=\delta_i-\delta_i'$. By Jensen's inequality,
$$
\E \exp \bigl(t \sum_{i=1}^n (\delta_i-\delta)a_i \bigr) \leq \E
\exp \bigl(t \sum_{i=1}^n (\delta_i-\delta_i')a_i \bigr) =
\prod_{i=1}^n \E \exp (t \tilde{\delta}_i a_i) = E.
$$

Set $\dd=\delta(1-\delta)$ and note that $\dd_i$ is $0$ with
probability $1 - 2 \dd$, and $1$ and $-1$, each with probability
$\dd$. Therefore,
$$
\E \exp (t \dd_i a_i) = (1 - 2 \dd) + \dd e^{t a_i} + \dd e^{- t
a_i} =  1 + 2 \dd ( \cosh (t a_i) - 1 ).
$$
Since $\cosh x  \le  1 + \frac{1}{2}x^2 e^{|x|}$ for all real $x$,
then
$$
E \le  \prod_{i=1}^n ( 1 + \dd t^2 a_i^2 e^{t |a_i|} )
  \le  \prod_{i=1}^n \exp( \dd t^2 a_i^2 e^{t |a_i|} )
  =    \exp \Big( \dd t^2 n \cdot \frac{1}{n} \sum_{i=1}^n a_i^2 e^{t
                 |a_i|} \Big).
$$
The normalized sum is estimated by Cauchy-Schwartz and using the
fact that $2 t \le 1$:
\begin{align*}
\frac{1}{n} \sum_{i=1}^n a_i^2 e^{t |a_i|}
  &\le  \Big( \frac{1}{n}  \sum_{i=1}^n |a_i|^4 \Big)^{\frac{1}{2}}
        \Big( \frac{1}{n}  \sum_{i=1}^n  e^{2 t |a_i|}
              \Big)^{\frac{1}{2}} \\
  &\le  \|a\|_{L_4^n}^2 \Big( \frac{1}{n} \; \sum_{i=1}^n  e^{|a_i|}
              \Big)^{\frac{1}{2}} \\
  &\le  2\|a\|_{L_4^n}^2
  \le  C,
\end{align*}
because $c \|a\|_{L_4^n}  \le  \|a\|_{\psi_1^n}  \le  1$. Hence,
$$
E  \le  \exp(C \dd t^2 n)  \le  \exp(C' \d t^2 n).
$$
We put this in a form convenient for applying Lemma
\ref{Chebychev}:
$$
\E \exp \Big( t \d n \cdot \frac{1}{\d n} \sum_{i=1}^n (\d_i-\d)
a_i \Big) \le  \exp \Big( \frac{C'}{\d n} (t \d n)^2 \Big)
$$
and apply the lemma for $\l = t \d n$ and $b^2 = \frac{C'}{\d n}$.
It follows that for every $t>0$,
$$
\P \Big\{ \frac{1}{\d n} \sum_{i=1}^n (\d_i-\d) a_i  >  2 c t \}
  \le  \exp( - c \d t^2 n),
$$
which completes the proof.
\endproof

\begin{Corollary} \label{remark}
Applying Proposition \ref{prop} to $-a_i$ it is evident that
$$
\P \Big\{ \Big| \sum_{i=1}^n (\d_i - \d) a_i \Big|
> t \d n\Big\}
  \le  2 \exp \Big( - \frac{c t^2 \d n}{M^2} \Big)
$$
where $M=\|a\|_{\psi_1^n}$ and $0 < t < M/2$.
\end{Corollary}

An easy application of this corollary is the fact that the
$\psi_2^n$-norm of points on the sphere determines the cardinality
of an almost isometric projection.

\begin{Corollary} \label{cor:projection-L_2}
There is an absolute constant $C$ for which the following holds.
For every integer $n$, any $f \in S(L_2^n)$ and every $\eps>0$, a
random set $\s \subset \{1,...,n\}$ of average cardinality
$(CM/\eps)^2$ satisfies with probability at least $1/2$ that
$$
1-\eps \leq \|P_\sigma f\|_{L_2^\s} \leq 1+\eps,
$$
where $M=\|f\|_{\psi_2^n}$.
\end{Corollary}

\proof The proof follows immediately from Corollary \ref{remark},
by taking $a_i=f^2(i)$ and $\delta n = (C M/\eps)^2$, and applying
\eqref{eq:app-Bern}.
\endproof

Note that a similar result can be easily derived for any $1 \leq p
< \infty$, simply by the fact that
$\|(a_i)\|_{\psi_p^n}=\|(a_i^p)\|_{\psi_1^n}$.

Corollary \ref{remark} can be used to present a new insight to the
well known Johnson-Lindenstrauss ``Flattening" Lemma \cite{JL},
which states that every set $\{x_1,...,x_n\} \subset \ell_2^n$ can
be $1+\eps$ isometrically embedded in $\ell_2^m$, where $m \leq
(C/\eps)^2 \log n$. One can formulate the Johnson-Lindenstrauss
Lemma as follows:

\begin{Theorem} \label{thm:J-L}
There is an absolute constant $C$ for which the following holds.
For every $f_1,...,f_n \in S(L_2^n)$ and every $\eps>0$ there is
an orthogonal operator $O$ and a set $\s \subset \{1,...,n\}$ of
cardinality at most $(C/\eps)^2 \log n$, such that for all $1 \leq
i \leq n$,
$$
1-\eps \leq \|P_{\s}Of_i\|_{L_2^\s} \leq 1+\eps.
$$
\end{Theorem}

As Corollary \ref{cor:projection-L_2} shows, an almost isometric
coordinate projection of $f$ is possible, as long as
$\|f\|_{\psi_2^n}$ is small; hence, the $\psi_2^n$ norm defines a
``good region" on the sphere for which a random coordinate
projection will be an almost isometry. In a similar way, this can
also be performed with many functions simultaneously:

\begin{Corollary}                             \label{cor:coord}
There is an absolute constant $C$ for which the following holds.
For every $f_1,...,f_n \in S(L_2^n)$ and every $\eps > 0$ a random
set $\s \subset \{1,...,n\}$ of cardinality $(CM/\eps)^2 \log n$
satisfies that with probability at least $1/2$,
$$
1-\eps \leq \|P_{\s}f_i\|_{L_2^\s} \leq 1+\eps, \ \ \ \ 1 \leq i
\leq n.
$$
where $M=\max_i\|f_i\|_{\psi_2^n}$.
\end{Corollary}

\proof As in Corollary \ref{cor:projection-L_2}, but taking
$\delta n = (C M/\eps)^2 \log n$, we obtain then for every $1 \le
i \le n$
$$
Pr \{1-\eps \leq \|P_{\s}f_i\|_{L_2^\s} \leq 1+\eps \} \ge 1 -
\frac{1}{2n}.
$$
Then
$$
Pr \{ \forall 1 \le i \le n, \ \ 1-\eps \leq
\|P_{\s}f_i\|_{L_2^\s} \leq 1+\eps \} \ge 1/2,
$$
which completes the proof.
\endproof

The connection to the Johnson-Lindenstrauss Lemma is easy: with
high probability, a random orthogonal operator $O$ will map any
set of $n$ vectors on the sphere to the ``good region", i.e. to
the region where the $\psi_2^n$ norm is bounded by an absolute
constant.

\begin{Lemma} \label{lemma:sphere}
There is an absolute constant $C$ such that for every integer $n$
and any $x \in S(L_2^n)$,
$$
Pr_{O_n} \bigl\{ \|Ox\|_{\psi_2^n} \geq C \bigr\} < \frac{1}{2n},
$$
where the probability measure is the Haar measure on the
orthogonal group.
\end{Lemma}

As a consequence, for every $f_1,...,f_n \in S(L_2^n)$,
$$
\max_i\|Of_i\|_{\psi_2^n} \leq C
$$
with probability greater than $1/2$, and thus Theorem
\ref{thm:J-L} is implied by Corollary \ref{cor:coord}.

\proof Clearly, it suffices to show that there is an absolute
constant $C$ such that
$$
Pr \bigl\{x \in S^{n-1}: \ \|x\|_{\psi_2^n} \geq
\frac{C}{\sqrt{n}} \bigr\} \leq \frac{1}{2n}.
$$
Consider the function $g:S^{n-1} \to \R$ defined by
$g(x)=\|x\|_{\psi_2^n}$. To estimate its Lipschitz constant,
observe that for every $x \in S^{n-1}$, $\|x\|_{\psi_2^n} \leq
\sqrt{2/\log{n}}$. Indeed, for $0 \leq x \leq 1$, $n^{x^2/2} \leq
nx^2+1$; hence,
$$
\frac{1}{n}\sum_{i=1}^n \exp\Bigl(\frac{x_i^2}{2}
\log{n}\Bigr)=\frac{1}{n}\sum_{i=1}^n n^{x_i^2/2} \leq
\frac{1}{n}\sum_{i=1}^n (nx_i^2+1) \leq 2.
$$

To bound the expectation of $g$ (with respect to the Haar measure
on the sphere), recall the median of the function
$f(x)=\sqrt{n}|x_1|$ satisfies that $M_f \sim c$, and that
$\|f\|_{{\rm lip}} \leq 1$. Hence, by concentration of measure on
the sphere \cite{MiS}, for any $s>c$ and every $1 \leq i \leq n$,
$$
Pr \bigl\{ x \in S^{n-1}: \sqrt{n}|x_i| \geq 2s \bigr\} \leq
\sqrt{\frac{\pi}{2}}e^{-s^2/2},
$$
and thus, $\E \exp(cnx_i^2) \leq 2$ for an appropriate absolute
constant $c$. Recall that there is an absolute constant $K$ such
that for every function $f$, $\|f\|_{\psi_2} \leq K \E \exp(f^2)$;
therefore, for $x=(x_1,...,x_n)$,
$$
\|\sqrt{cn}x\|_{\psi_2^n} \leq \frac{K}{n}\sum_{i=1}^n
\exp(cnx_i^2).
$$
Taking the expectation with respect to $x$ on the sphere,
$$
\E \|\sqrt{n}x\|_{\psi_2^n} \leq \frac{K}{n}\sum_{i=1}^n \E
\exp(cnx_i^2) \leq K'
$$
for an absolute constant $K'$.

By the concentration of measure on the sphere applied to the
function $g$,
$$
Pr \bigl\{x \in S^{n-1} : \ \|x\|_{\psi_2^n} \geq
\frac{C}{\sqrt{n}}+t\bigr\} \leq \sqrt{\frac{\pi}{2}} e^{-ct^2n
\log{n}},
$$
and the claim follows by selecting $t=C'/\sqrt{n}$.
\endproof

\end{document}